\documentclass[12pt]{amsart}

\usepackage[textheight=240mm, textwidth=169mm, top=20mm, centering]{geometry}

\usepackage{amsmath, amssymb, amsthm, amscd, mathrsfs, mathtools}

\usepackage{graphicx, hhline}
\usepackage[all]{xy}

\usepackage{enumerate}
\usepackage{xcolor}
\usepackage{url}
\usepackage{bm}

\numberwithin{equation}{section}

\theoremstyle{plain}
\newtheorem{thm}{Theorem}[section]
\newtheorem{prop}[thm]{Proposition}
\newtheorem{lem}[thm]{Lemma}

\newtheorem{cor}[thm]{Corollary}

\theoremstyle{definition}
\newtheorem{dfn}[thm]{Definition}

\newtheorem{notation}[thm]{Notation}

\newtheorem{rem}[thm]{Remark}

\begin{document}

\title[Nef Direct Image Sheaf]
{Metric nefness of direct image sheaves under \\
weakly semi-stable K\"ahler morphisms}

\author{Osamu Fujino}
\author{Shigeharu Takayama}

\address{O.~F. Department of Mathematics, Graduate School of Science, Kyoto University, Kyoto 606-8502, Japan.}
\email{fujino@math.kyoto-u.ac.jp}

\address{S.~T. Graduate School of Mathematical Sciences, the University of Tokyo, 3-8-1 Komaba, Meguro-ku, Tokyo 153-8914, Japan.}
\email{taka@ms.u-tokyo.ac.jp}

\begin{abstract}
We investigate the nefness of the direct image sheaves of relative adjoint bundles under weakly semi-stable K\"ahler morphisms between complex spaces.
We show that the direct image sheaves are locally free and that their  tautological line bundles are metrically nef, strengthening the previously known results on torsion-freeness and certain weak positivity.
\end{abstract}

\date{\today}
\makeatletter 
\@namedef{subjclassname@2020}{\textup{2020} Mathematics Subject Classification}



\date{\today}
\maketitle

\baselineskip=16pt

\linespread{1.1} 

\makeatletter
\@namedef{subjclassname@2020}{\textup{2020} Mathematics Subject Classification}
\makeatother

\newcommand{\disp}{\displaystyle}

\newcommand{\BC}{{\mathbb{C}}}
\newcommand{\BN}{{\mathbb{N}}}
\newcommand{\BP}{{\mathbb{P}}} 
\newcommand{\BQ}{{\mathbb{Q}}}
\newcommand{\BR}{{\mathbb{R}}}
\newcommand{\BZ}{{\mathbb{Z}}}

\newcommand{\CE}{{\mathcal{E}}}
\newcommand{\CF}{{\mathcal{F}}}
\newcommand{\CI}{{\mathcal{I}}}
\newcommand{\CJ}{{\mathcal{J}}}
\newcommand{\CL}{{\mathcal{L}}}
\newcommand{\CM}{{\mathcal{M}}}
\newcommand{\CO}{{\mathcal{O}}}

\renewcommand{\OE}[1]{{\mathcal{O}_{\mathbb{P}(E)}(#1)}}

\newcommand{\bfx}{\mathbf {x}}
\newcommand{\bfy}{\mathbf {y}}
\newcommand{\bfz}{\mathbf {z}}
\newcommand{\bfo}{\mathbf {o}} 

\newcommand{\ora}{\overrightarrow}

\newcommand{\bfgc}{\mathbf {\gc}}
\newcommand{\bflam}{\mathbf {\lambda}}
\newcommand{\bfth}{\mathbf {\theta}}
\newcommand{\fm}{{\mathfrak {m}}}

\newcommand{\ga}{{\alpha}}
\newcommand{\gb}{{\beta}}
\newcommand{\gc}{{\gamma}}
\newcommand\ep{\varepsilon}
\newcommand{\Ga}{{\Gamma}}
\newcommand{\Del}{{\Delta}}
\newcommand{\Dl}{{\Delta}}
\newcommand{\del}{{\delta}}

\newcommand{\lam}{{\lambda}}
\newcommand{\Lam}{{\Lambda}}
\newcommand{\vph}{{\varphi}}
\newcommand{\sg}{{\sigma}}
\newcommand{\Th}{{\Theta}}
\newcommand\w{\omega}
\newcommand\Om{\Omega}

\newcommand\wtil{\widetilde}
\newcommand\what{\widehat}

\newcommand\tU{\widetilde U}
\newcommand\tx{\widetilde x}
\newcommand\tX{\widetilde X}

\newcommand\sm{\setminus}
\newcommand\ol{\overline}
\newcommand\ovl{\overline}
\newcommand\ot{\otimes}
\newcommand\wed{\wedge}

\newcommand\lra{\longrightarrow}

\newcommand\ai{\sqrt{-1}}
\newcommand\rd{{\partial}}
\newcommand\rdb{{\overline{\partial}}}

\newcommand\cpt{{\rm cpt}}
\newcommand\Supp{\mbox{{\rm Supp}}\, }

\newcommand{\qedd}{\hfill {\Large $\triangleleft$}} 


\section{Introduction}

We are interested in various positivity properties of direct image sheaves of adjoint bundles.
In particular, in this paper, we study the case where direct image sheaves become nef vector bundles.
Our main result is stated as follows.

\begin{thm}\label{mt}
Let $X$ be a complex space, $Y$ a complex manifold, and $f:X\to Y$ a proper, surjective, locally K\"ahler holomorphic map with connected fibers.
Suppose that $f$ is weakly semi-stable and that either the toroidal structures of $X$ and $Y$ are strict or $X$ has at worst quotient singularities.
Let $M$ be a holomorphic vector bundle on $X$ equipped with a Hermitian metric $h$ on $X$ that is Nakano semi-positive on the regular locus $X_{reg}$.
Then, the direct image sheaf $E:=f_{*}(K_{X/Y} \otimes M)$ is locally free, and the canonical $L^{2}$-metric $g$ on $E$ is singularly Griffiths semi-positive and has at worst logarithmic divergence.
In particular, the universal quotient line bundle $\OE1$ on $\BP(E)$ is metrically nef.
\end{thm}

It is known, to some extent, in the setting of \ref{mt} that $E=f_{*}(K_{X/Y} \otimes M)$ is reflexive and singularly Griffiths semi-positive with respect to the canonical $L^2$-metric $g$ which is given by fiberwise integrals.
The weak semi-stability will be used to prove both that $E$ is locally free and that the singularities of $g$ are mild enough to conclude the metric nefness.
The probably unfamiliar notion of metric nefness is defined as follows.

\begin{dfn}[{cf.~\cite[(6.10)]{DLN}}]\label{mnef}
Let $X$ be a complex manifold and $\omega$ a reference Hermitian metric on $X$.
A holomorphic line bundle $L$ on $X$ is said to be {\it metrically nef} if, for any relatively compact open subset $W \subset X$ and any constant $\varepsilon > 0$, there exists a Hermitian metric $h_{W}$ on $L|_W$ such that the curvature form satisfies $\sqrt{-1}\Theta_{h_{W}} > -\varepsilon\omega$ on $W$.
\qedd
\end{dfn}

Note that $X$ is not necessarily compact or algebraic, and that metric nefness and nefness in algebraic geometry coincide on smooth projective varieties.
When $X$ and $Y$ are projective and $M$ is of rank $1$, \ref{mt} was essentially proved by the first-named author in \cite[Remarks 4.3.9, 4.3.10, and 4.3.13]{fujino3}.
His method is based on Viehweg's arguments, namely, covering (base-change) tricks and vanishing theorems.
These are algebraic arguments and are very different from our approach in this paper.
One major motivation for \ref{mt} comes from recent developments in the study of projective morphisms between complex spaces and the minimal model program in the K\"ahler category (see, e.g., \cite{fujino4}, \cite{hashizume}).
We now state a corollary using a formulation that is familiar to algebraic geometers.

\begin{cor}
Let $f\colon X\to Y$ be a proper surjective morphism of smooth projective varieties with connected fibers, and let $\mathcal{L}$ be a semi-positive line bundle on $X$. 
By the weakly semi-stable reduction theorem, we can construct a commutative diagram 
$$
\xymatrix{
X' \ar[r]^\rho \ar[d]_{f'} & X \ar[d]^f \\ 
Y' \ar[r]_\tau & Y,
}
$$ 
where $\tau\colon Y'\to Y$ is a generically finite surjective morphism of smooth projective varieties, and $f'\colon X'\to Y'$ is a weakly semi-stable morphism between projective varieties such that $X'$ is a partial resolution of the dominant components of $X\times_Y Y'$ and has only quotient singularities. 
Set $\mathcal{L}':=\rho^*\mathcal{L}$. 
Then 
$$
\mathcal{E} := f'_*\left(K_{X'/Y'}\otimes \mathcal{L}'\right)
$$ 
is a nef locally free sheaf on $Y'$. 

More precisely, let $\pi\colon \mathbb{P}_{Y'}(\mathcal{E})\to Y'$ be the natural projection. 
The morphism $f'$ is toroidal with respect to toroidal structures $(U_{X'} \subset X')$ and $(U_{Y'}\subset Y')$.
Then, the line bundle $\mathcal{O}_{\mathbb{P}_{Y'}(\mathcal{E})}(1)$ admits a singular Hermitian metric $h$ with semi-positive curvature current such that $h$ is smooth on $\pi^{-1}(U_{Y'})$, and the Lelong number (see \ref{LelongNo} for the definition) of the local weight of $h$ is zero everywhere. 
\end{cor}

We provide further technical comments on \ref{mt}.
We will prove that all direct image sheaves $R^{i}f_{*}(K_{X/Y} \otimes M)$ are locally free. 
They are already known to be torsion-free \cite{takegoshi}, which implies local freeness and the constancy of fiberwise dimensions $h^{i}(X_{y}, (K_{X/Y} \otimes M)|_{X_y})$ when $\dim Y = 1$ by the base change theorem. 
The general case ($\dim Y > 1$) reduces to this 1-dimensional case. 
Indeed, the weak semi-stability is preserved by restriction to a general smooth curve $C \subset Y$, meaning that $f_{C}: X_{C} := X \times_{Y} C \to C$ remains weakly semi-stable. 
Applying the base change argument over varying $C$ implies that the dimensions $h^{i}(X_{y}, (K_{X/Y} \otimes M)|_{X_y})$ are constant for all $y \in Y$.

The canonical $L^2$-metric $g$ on $E$ canonically induces a singular Hermitian metric $g_{\OE1}$ with a semi-positive curvature current on the line bundle $\OE1$ over the projective space bundle $\pi : \BP(E) \to Y$. 
Weak semi-stability provides appropriate local coordinates for calculating and estimating $g$, showing that the singularities of both $g$ and $g_{\OE1}$ are sufficiently mild. 
We then deduce the metric nefness of $\OE1$ from Demailly's regularization theorem \cite{D92}, which is recalled below.

Let us discuss metric nefness in general.
A Cartier divisor $D$ on a proper algebraic variety is said to be \textit{nef} if the intersection number $D \cdot C$ is semi-positive for any curve $C$ ($D \cdot C \ge 0$).
An analytic counterpart of nefness can be formulated (see, e.g., \cite[(6.10)]{DLN}).
Definition \ref{mnef} is a modified version of it, adapted to general complex manifolds.
The property of being metrically nef in \ref{mnef} does not depend on the choice of the reference Hermitian metric $\omega$ there.
If a line bundle $L$ is metrically nef and if $X$ contains an irreducible compact complex subspace $V$ of dimension $d>0$, then the intersection number is semi-positive, i.e., $L^{d} \cdot V \ge 0$.
This means that $L$ is nef in the sense of algebraic geometry if $X$ is a smooth proper algebraic variety.
It is clear that if a line bundle $L$ admits a Hermitian metric with semi-positive curvature, then $L$ is metrically nef.
However, the converse does not hold in general (see, e.g., \cite[Example 1.7]{DPS}, \cite[\S 5]{fujino0}).
One would need, in principle, a smooth Hermitian metric on $L$ to test the nefness condition $L \cdot C \ge 0$.
On the other hand, canonical metrics on vector bundles (such as the canonical $L^2$-metric $g$ on $E$ in \ref{mt}) can be singular in general.
To address this issue, Demailly relaxed the smoothness assumption on $h$ to conclude the metric nefness of $L$.

\begin{prop}[{cf.~\cite[Corollary 6.4]{D92}}]\label{Dcor64}
Suppose that a holomorphic line bundle $L$ on a complex manifold $X$ admits a singular Hermitian metric $h$ with semi-positive curvature, whose local plurisubharmonic  potential functions have vanishing Lelong numbers everywhere on $X$ except possibly at a countable set of points.
Then $L$ is metrically nef.
\end{prop}

Returning to the situation of \ref{mt}, we shall show that the line bundle $\OE1$ equipped with $g_{\OE1}$ satisfies the condition in \ref{Dcor64}.
Consequently, our primary task is to bound the singularities of the metric $g$ on $E$.
Although only a very rough classification of the potential singularities of canonical $L^{2}$-metrics is known, we introduce a ``logarithmic divergence property'' to characterize the mildness of the singularities of singular Hermitian metrics.

The organization of this paper is as follows.
In \S 2, we explain weakly semi-stable morphisms and prove the local freeness of all direct images $R^{i}f_{*}(K_{X/Y}\ot M)$.
After some preliminaries, including the definition of Lelong numbers, we give an outline of the proof of \ref{Dcor64} in \S 3.
In \S 4, we discuss singular Hermitian metrics with logarithmic singularities in a general setting and formulate a criterion for the line bundle $\OE1$ to be metrically nef.
We then prove in \S 5 that canonical $L^{2}$-metrics possess such logarithmic singularities.
In \S 6, we give a few remarks on Nakano semi-positivity on possibly singular spaces.

In this paper, we define a holomorphic map $f:X \to Y$ to be {\it locally K\"ahler} if for every point $y \in Y$, there exists an open neighborhood $V$ of $y$ such that its preimage $f^{-1}(V)$ is K\"ahler (cf.~\cite[6.1]{takegoshi}).
See \ref{Nsp} for the definition of Nakano semi-positivity. 
Any alternative definition yields the same result; see \ref{any-resol}--\ref{any-pb}.

Acknowledgments.
The first author was partially supported by JSPS KAKENHI Grant Number JP23K20787.
The second author was supported by JSPS KAKENHI Grant Number JP26H01995.

AI-based tools were used for English language editing and proofreading of the manuscript.

\section{Weak semi-stability and local freeness}

We recall the definition of weakly semi-stable morphisms and prove the local freeness of all direct image sheaves $R^i f_*(K_{X/Y}\otimes M )$ in \ref{mt}.

\begin{dfn}[{Weakly semi-stable morphisms, see \cite{ak} and \cite{eh}}]\label{def-wss}
Let $f\colon X\to Y$ be a proper morphism of complex spaces. 
We say that $f$ is \emph{weakly semi-stable} if 
\begin{itemize}
\item[(i)] $X$ and $Y$ admit toroidal structures $(U_X\subset X)$ and $(U_Y\subset Y)$ with $f^{-1}(U_Y)=U_X$, 
\item[(ii)] $f$ is toroidal with respect to these toroidal structures, 
\item[(iii)] $f$ is equi-dimensional, 
\item[(iv)] all the fibers of $f$ are reduced, and 
\item[(v)] $Y$ is smooth.  
\end{itemize}
\end{dfn}

We do not need to suppose that $f$ has connected fibers.
We do not recall further technical details concerning toroidal embeddings (see \cite{ak}, \cite{cls}).
By a basic property of weak semi-stability, $X$ has only rational Gorenstein singularities (see \cite[Proposition 6.4]{ak}), which in turn implies that $X$ has only canonical Gorenstein singularities.
Here are a couple of further remarks.

\begin{rem}
In \ref{def-wss}, 
the condition $f^{-1}(U_Y)=U_X$ implies that $f$ is smooth over $U_Y$ because 
$f$ is a toroidal morphism with respect to the toroidal structures 
$(U_X\subset X)$ and $(U_Y\subset Y)$. 
\qedd
\end{rem}

\begin{rem}
Let $(U_X\subset X)$ be a toroidal embedding. If every irreducible component of the boundary $X\setminus U_X$ is normal, then $(U_X\subset X)$ is called a \emph{toroidal embedding without self-intersection} or a \emph{strict toroidal embedding}. 

Working with strict toroidal embeddings is necessary to properly define the associated polyhedral complexes, and it allows us to avoid inessential technical complications.
\qedd
\end{rem}

In \ref{def-wss}, we do not assume that $X$ has at worst quotient singularities. 
However, such a quasi-smooth model can be achieved via a modification, provided that the toroidal structures of $X$ and $Y$ are strict. 
Here, a {\it modification} means a proper, surjective, and generically one-to-one holomorphic map.

\begin{lem}\label{qs-wss}
Let $f:X\to Y$ be a weakly semi-stable and locally K\"ahler morphism between complex spaces, and suppose that the toroidal structures of $X$ and $Y$ are strict.
After restricting $f$ over any relatively compact open subset of $Y$, there exists a modification $\vph : X' \to X$ such that 
\begin{itemize}
\item[(i)] 
$X'$ is quasi-smooth, i.e., it has quotient singularities at worst, 
\item[(ii)] 
the morphism $\vph$ is projective and toroidal, 
\item[(iii)] 
the composition $f' = f \circ \vph :X'\to Y$ is identical to $f:X\to Y$ over $U_{Y}$, and
\item[(iv)] 
$f' : X'\to Y$ is weakly semi-stable and locally K\"ahler.
Moreover, the toroidal structure of $X'$ is strict.
\end{itemize}
\end{lem}

Although the following proof is more or less well known to experts, we include it for the sake of completeness. 

\begin{proof}
Throughout the proof, we repeatedly replace $Y$ with a relatively compact open subset of $Y$. 

By assumption, $(U_X\subset X)$ and $(U_Y\subset Y)$ are strict toroidal embeddings and $f\colon X\to Y$ is toroidal with respect to these toroidal structures. 
We can consider rational conical polyhedral complexes $\Delta_X$ and $\Delta_Y$ associated to $(U_X\subset X)$ and $(U_Y\subset Y)$, respectively. 
For details, see \cite[\S 1.3]{ak} and \cite[\S 2.1]{eh}.
We may assume that these complexes $\Delta_X$ and $\Delta_Y$ are finite by replacing $Y$ with a relatively compact open subset.  
In this setting, we have a morphism of rational conical polyhedral complexes $f_\Delta\colon \Delta_X\to \Delta_Y$ associated to $f\colon X\to Y$. 

By \cite[Remark 4.5]{ak}, a finite sequence of suitable star subdivisions of $\Delta_X$ yields a simplicial subdivision $\Delta_{X'}$ without adding any new $1$-dimensional cones. For details, see, for example, \cite[Proposition 11.1.7]{cls}.
Let $\varphi\colon X'\to X$ be the bimeromorphic toroidal morphism corresponding to this subdivision. 
Since $\varphi$ is obtained by a finite sequence of star subdivisions, it becomes projective by shrinking $Y$ to a relatively compact open subset. 
Furthermore, because $\Delta_{X'}$ is simplicial, $X'$ has only quotient singularities. 

Since $f$ is weakly semi-stable, it is equi-dimensional, which implies that $f_\Delta\colon \Delta_X\to \Delta_Y$ maps cones onto cones (see \cite[Lemma 4.1]{ak}). 
Because the $1$-dimensional cones of $\Delta_X$ and $\Delta_{X'}$ are exactly the same, 
the induced morphism $\Delta_{X'}\to \Delta_Y$ also maps cones onto cones. 
Hence, the induced map $f'\colon X'\to Y$ is equi-dimensional (see \cite[Lemma 4.1 and Remark 4.2]{ak}). 

Since $Y$ is smooth, the equi-dimensionality of $f'$ implies that $f'$ is flat. 
Finally, by \cite[Lemma 5.2 and Remark 5.4]{ak}, this flatness combined with the fact that $f\colon X\to Y$ is weakly semi-stable guarantees that $f'\colon X'\to Y$ is also weakly semi-stable. 
As both $f$ and $\vph$ are locally K\"ahler, the composition $f'=f\circ \vph$ is locally K\"ahler (\cite[Proposition 6.2]{takegoshi}).
\end{proof}

The following remark is clear from the proof of \ref{qs-wss}.

\begin{rem}\label{rem-quotient}
In the setting of \ref{qs-wss}, if $f\colon X\to Y$ is algebraic, 
we can construct a projective birational toroidal morphism $\varphi\colon X'\to X$ without replacing $Y$ with a relatively compact open subset. 
\end{rem}

We now proceed to discuss the local freeness of direct image sheaves. 
Note that we assume neither that the toroidal structures of $X$ and $Y$ are strict nor that $X$ has at worst quotient singularities.

\begin{lem}[Torsion-freeness]\label{lem-torsion-free}
Let $f\colon X\to Y$ be a weakly semi-stable locally K\"ahler morphism, and $M$ a Nakano semi-positive vector bundle on $X$. 
Then $R^i f_*(K_{X/Y}\otimes M)$ is torsion-free for every $i \ge 0$.
\end{lem}

\begin{proof}
Since the assertion is local on $Y$, we may freely replace $Y$ with a relatively compact open subset throughout the proof without explicit mention. 
In particular, we may assume that $K_Y \simeq \mathcal{O}_Y$, which yields $K_{X/Y} \simeq K_X$. 
We also note that $X$ has only canonical Gorenstein singularities. 
Thanks to the Grauert--Riemenschneider (relative Kawamata--Viehweg) vanishing theorem, we may replace $(X,M)$ with $(\wtil X,\mu^{*}M)$, where $\mu:\wtil X\to X$ is a resolution of singularities. 
Thus, we may further assume that $X$ is a smooth K\"ahler manifold. 
In this smooth setting, it is already known by \cite[Theorem 6.5]{takegoshi} that $R^i f_*(K_X \otimes M)$ is torsion-free for every $i$. 
\end{proof}

\begin{thm}[{Local freeness, see \cite[Remark 4.3.10]{fujino3}}]\label{main-thm}
Let $f\colon X\to Y$ be a weakly semi-stable locally K\"ahler morphism, and $M$ a Nakano semi-positive vector bundle on $X$. 
Then $R^i f_*(K_{X/Y}\otimes M)$ is locally free for every $i \ge 0$. 
In particular, $f_* (K_{X/Y}\otimes M)$ is locally free.  
\end{thm}

\begin{proof}
Fix an arbitrary point $P\in Y$. It suffices to show that $R^i f_*(K_{X/Y}\otimes M)$ is locally free in a neighborhood of $P$. 
Therefore, we may assume that $Y$ is a small open polydisc in $\mathbb{C}^n$ centered at $P$. 

Since $f$ is weakly semi-stable, it is flat. Moreover, $X$ has only Gorenstein singularities, so the relative canonical sheaf $K_{X/Y}$ is invertible. This implies that $\mathcal F:=K_{X/Y}\otimes M$ is locally free on $X$. Because $f$ is a flat morphism and $\mathcal{F}$ is locally free, $\mathcal{F}$ is flat over $Y$. 

We choose $n-1$ general smooth hyperplanes $H_1, \dots, H_{n-1}$ passing through $P$. By choosing $H_j$ in general position with respect to the stratification defining the weakly semi-stable morphism $f$, the complete intersection $C = H_1 \cap \dots \cap H_{n-1}$ is a smooth curve passing through $P$. Moreover, the restricted morphism 
$$
f_C \colon X_C \to C, \quad \text{where } X_C = X \times_Y C,
$$
remains a weakly semi-stable morphism by \cite[Lemma 6.2]{ak} (see also 
\cite[Lemma 2.12 (2)]{karu}). 

We put $\mathcal{F}_C := \mathcal{F}|_{X_C}$. Then we have 
$\mathcal F_C\cong K_{X_C/C} \otimes M|_{X_C}$ by adjunction. By Lemma \ref{lem-torsion-free}, the higher direct image $R^j(f_C)_*\mathcal{F}_C$ is locally free on $C$ for every $j \ge 0$. 
Since $\mathcal{F}_C$ is flat over $C$ and $R^j(f_C)_*\mathcal{F}_C$ is locally free for every $j \ge 0$, we can deduce that the base change morphism 
$$
R^j(f_C)_*\mathcal{F}_C \otimes \mathbb{C}(y) \to H^j(X_y, \mathcal{F}_y)
$$
is an isomorphism for all $y \in C$ and for all $j \ge 0$. This follows by applying \cite[Chapter III, Theorem 3.4]{bs} using descending induction on $j$, starting from sufficiently large $j$ where the cohomology groups naturally vanish. 
This implies that the function $h^i(y) := \dim_{\mathbb C} H^i(X_y, \mathcal{F}_y)$ is constant on $C$. 

Since $C$ is a general curve passing through $P$, it intersects the generic locus of $Y$. Thus, the constancy of $h^i(y)$ on $C$ implies that the value $h^i(P)$ equals the generic value of $h^i(y)$ in a neighborhood of $P$ in $Y$. 
By Grauert's semicontinuity theorem (see, for example, 
\cite[Chapter III, Theorem 4.12.(i)]{bs}), the function $y \mapsto h^i(y)$ is upper semicontinuous on $Y$. Thus, $h^i(y)$ must be constant in a sufficiently small open neighborhood of $P$ in $Y$. 

Finally, since $\mathcal{F}$ is flat over $Y$ and the function $h^i$ is locally constant on $Y$, Grauert's continuity theorem (see, for example, \cite[Chapter III, Theorem 4.12.(ii)]{bs}) ensures that $R^if_*\mathcal{F} = R^if_*(K_{X/Y}\otimes M)$ is locally free around $P$. This completes the proof.
\end{proof}

\section{A few analytic preliminaries}

We review a few basic properties of singular Hermitian metrics.
General references include \cite{DLN}, \cite{PT}, and \cite{HPS}.
Let $X$ be a complex manifold, and $E$ a holomorphic vector bundle on $X$ of rank $r$.
We refer the reader to \cite[\S 17]{HPS} for the definition of singular Hermitian metrics on vector bundles; here, we only recall the notions of curvature negativity and positivity.

\begin{dfn}
A singular Hermitian metric $h$ on $E$ is said to be {\it Griffiths semi-negative} if, for any local section $s \in H^{0}(U,E)$, the function $\log |s|_{h}$ is plurisubharmonic (psh for short) on $U$.
A singular Hermitian metric $h$ on $E$ is said to be {\it Griffiths semi-positive} if its dual $(E^{*},h^{*})$ is Griffiths semi-negative.
\qedd
\end{dfn}

There is a canonical way to induce a singular Hermitian metric on $\OE1$ from $(E, h)$.

\begin{rem}
Let $\pi : \mathbb{P}(E) \to X$ be the projective space bundle, and $\OE1$ the universal quotient line bundle.
We employ the following procedure to study the singularities of the induced metric on $\OE1$.

Suppose that a singular Hermitian metric $h$ on $E$ is Griffiths semi-positive.
Then its dual $h^{*}$ is a singular Hermitian metric on $E^{*}$ with Griffiths semi-negative curvature.
The pull-back $\pi^{*}E^{*}$ inherits a singular Hermitian metric $\pi^{*}h^{*}$ that also has Griffiths semi-negative curvature.
Consequently, the sub-bundle $\OE{-1} \hookrightarrow \pi^{*}E^{*}$ is equipped with a singular Hermitian metric with Griffiths semi-negative curvature by restriction:
$$
	h_{\OE{-1}} := \big( \pi^{*}h^{*} \big)|_{\OE{-1}}.
$$
Finally, we obtain a singular Hermitian metric $h_{\OE1} = {h_{\OE{-1}}}^{-1}$ on $\OE1$ with Griffiths semi-positive curvature.
\qedd
\end{rem}

The Lelong number is a quantity used to measure the analytic singularities of a psh function.

\begin{dfn}\label{LelongNo}
Let $\varphi$ be a psh function on an open set $U$ in $\mathbb{C}^n$. 
The {\textit{Lelong number}} $\nu(\varphi, x)$ of $\varphi$ at $x \in U$ is defined by 
$$
	\nu(\varphi, x)
	= \liminf_{z \to x} \frac{\varphi(z)}{ \log |z-x|} \, \in [0,+\infty).
$$
The Lelong number is independent of the choice of local coordinates and can, therefore, be defined on complex manifolds.
\qedd
\end{dfn}

The Lelong number is normalized so that $\nu(\log|z|,0)=1$.
A typical example with $\nu(\vph,0)=0$ is $\vph = - \log(-\log|z|)$.
This type of singularity will play a crucial role in our main results.


A possibly singular Hermitian metric $h$ on a line bundle is locally written as $h= e^{-\vph}$ with $\vph \in L^{1}_{loc}(U, \BR)$ on each small open set $U$.
The function $\vph$ is called a {\it local potential} of $h$.
The curvature current of $h$ is semi-positive on $U$ if and only if its local weight $\vph$ admits a psh representative.
We always choose this representative.
In this case, we can discuss the Lelong number $\nu(\vph,x)$ of $\vph$ at $x \in U$, which is associated with $h$ and independent of the choice of local trivializations.
The curvature current $\ai\Theta_{h}$ of $h$ is given by $\ai\rd\rdb \vph$ on $U$, where $\rd$ and $\rdb$ are taken in the sense of distributions.

 A locally integrable function $\vph \in L^{1}_{loc}(X,\BR)$ is said to be {\it quasi-psh} if at each $x\in X$, there exists a local chart $(U,z)$ around $x$ such that $\vph+\del_{U}|z|^{2}$ is psh on $U$ for some constant $\del_{U}>0$.
One can also define the Lelong numbers for quasi-psh functions.

We now give an outline of the proof of \ref{Dcor64}.
Namely, this serves as a guide on how to read \cite[6.4]{D92} for a relatively compact open set.
The original statement in \cite[6.4]{D92} is formulated for a compact complex manifold $X$.
However, it can be extended to general complex manifolds, as indicated at the top of \cite[p.~363]{D92}.
(Refer to \cite[Proposition 3.5]{FF} for a proof of \ref{Dcor64} based on Nadel's vanishing and Mumford's regularity when the variety $X$ is smooth projective.)

\begin{proof}[Outline of the proof of \ref{Dcor64}]
We fix a Hermitian form $\w$ on $X$.
Let $\del >0$ be an arbitrarily small constant.
We choose relatively compact open sets $W \subset V \subset X$ such that $\ol W \subset V$.

Let $g$ be a Hermitian metric on $L$.
Then, we can write $h = ge^{-\vph}$ for some $\vph \in L^{1}_{loc}(X, \BR)$.
Thus, we have $\ai\rd\rdb \vph = \ai \Theta_{h} - \ai\Theta_{g} \ge 0+\gc$, where $\gc = - \ai\Theta_{g}$ is a real $(1,1)$-form on $X$.
Consequently, we can suppose $\vph$ to be quasi-psh on $X$.
For any constant $c>0$, we set $E_{c}(\vph) = \{x \in X \mid \nu(\vph,x)\ge c\}$, which is the upper-level set of the Lelong number.
By Siu's analyticity theorem \cite{S}, $E_{c}(\vph)$ is an analytic set in $X$.
In particular, by our assumption on the Lelong number, $E_{c}(\vph) \cap V$ consists of finitely many points in $V$.
We then proceed to regularize $\vph$ with a small loss of positivity.

We now apply \cite[3.7]{D92} and \cite[(13.12)]{DLN} (cf.~\cite[1.1]{D92}) to the relatively compact set $V$ as indicated at the top of \cite[p.~363]{D92}.
This yields a sequence of quasi-psh functions $\vph_{c,k}$ ($k=1,2,\ldots$) on $V$ satisfying a number of conditions.
First of all, each $\vph_{c,k}$ is smooth on $V \setminus E_{c}(\vph)$ (and has the so-called analytic singularities).
Furthermore, we have $\ai\rd\rdb \vph_{c,k} \ge \gc - \ep_{k}\w$ on $V$ for some decreasing sequence of positive numbers $\{\ep_{k}\}_{k}$ such that $\lim \ep_{k}=0$.
We set $h_{c,k}= g e^{-\vph_{c,k}}$, which defines a singular Hermitian metric on $L|_{V}$.
We then fix a sufficiently large integer $k$ so that $0 < \ep_{k} < \del$ holds.
In particular, we obtain $\ai\Theta_{h_{c,k}} > - \del\w$ on $V$.

We consider $ \vph_{V} = \vph_{c,k}$ on $V$, which is smooth outside the finite set $E_{c}(\vph) \cap V$ and satisfies $\ai\rd\rdb  \vph_{V} > \gc - \del \w$ on $V$.
Then, \cite[6.3]{D92} (using a classical smoothing procedure) shows that there exists a smooth function $\vph_W$ on $W$ such that $\ai\rd\rdb  \vph_{W} > \ai\rd\rdb  \vph_{V} -\del \w > \gc - 2\del\w$ on $W$.
Here, we have passed to a relatively compact subset $W \subset V$ and lost some positivity again.
Consequently, the smooth Hermitian metric $h_{W} = ge^{- \vph_{W}}$ on $L|_{W}$ satisfies $\ai \Theta_{h_{W}}= \ai\Theta_{g} + \ai\rd\rdb  \vph_{W} > -2\del \w$ on $W$.
\end{proof}


\section{Singular Hermitian metrics with logarithmic singularities}

We introduce a specific class of singular Hermitian metrics, namely those with logarithmic singularities at worst.
Here, logarithmic singularities refer to logarithmic growth of the metric coefficients, rather than to logarithmic poles of their local weights.
This type of singularity appears in various other contexts.
We fix the following notation for this section.

\begin{notation}\label{adm}
Let $X$ be a complex manifold of dimension $n$, and $D$ a reduced normal crossing divisor on $X$.
For each $x \in D$, we can choose a local coordinate system $(U, t=(t_{1},\ldots, t_{n}))$ on $X$ centered at $x$ such that $D \cap U = \{\prod_{j=1}^{n_{x}} t_{j}=0\}$ for some integer $n_{x}>0$.
We identify $(U,t)$ with a subset of $\BC^{n}$ and set $|t| = (\sum |t_{i}|^{2})^{1/2}$.

Let $E$ be a holomorphic vector bundle on $X$ of rank $r$ equipped with a singular Hermitian metric $h$ which is continuous and positive definite on $X \sm D$.
We also assume that $E$ is trivial on the local chart $U$ chosen above, and choose a local frame $e_{1},\ldots,e_{r} \in H^{0}(U,E)$ such that $E|_{U}=\bigoplus_{i=1}^{r} \CO_{U} e_{i}$.
We express $h$ as a matrix-valued function on $U$ by
$$
	h = ( h_{i\ol j})_{1 \le i, j \le r}, \quad h_{i\ol j} = h(e_{i},e_{j}).
$$
Let $0< \lam_{1}(t) \le \lam_{2}(t) \le \ldots \le \lam_{r}(t) <+\infty$ denote the eigenvalues of $(h_{i\ol j})$ at $t \in U \setminus D$.

In this setting, we call $(U, t)$ an ``admissible chart''.
\qedd
\end{notation}

As our discussion in this section is local, the normal crossing divisor $D$ need not be a simple normal crossing divisor.
We will discuss the singularities of $h$ along $D$.

\begin{dfn}\label{LDD}
Let $(X,D)$ and $(E,h)$ be as in \ref{adm}, and let $k$ be a non-negative integer.

(1)
A singular Hermitian metric $h$ on $E$, or the pair $(E,h)$, is said to {\it have a logarithmic divergence of order at most $k$ along $D$}
if, for every $x \in D$ and every admissible chart $(U,t=(t_{1},\ldots,t_{n}))$ centered at $x$, 
there exist constants $c_{x}>0$ and $\delta_{x} \in (0,1)$ such that, for any $t \in \{|t| < \delta_{x}\} \sm D$,
$$ 
	\lam_{r}(t) \le c_{x} \prod_{j=1}^{n_{x}} (- \log|t_{j}|)^{k} 
$$ 
holds.

(2) 
A singular Hermitian metric $h$ on $E$, or the pair $(E,h)$, is said to {\it have a logarithmic decay of order at most $-k$ along $D$}
if, for every $x \in D$ and every admissible chart $(U,t=(t_{1},\ldots,t_{n}))$ centered at $x$, 
there exist constants $c_{x}>0$ and $\delta_{x} \in (0,1)$ such that, for any $t \in \{|t| < \delta_{x}\} \sm D$,
$$ 
	\lam_{1}(t) \ge c_{x} \prod_{j=1}^{n_{x}} (- \log|t_{j}|)^{-k} 
$$ 
holds.
\qedd
\end{dfn}

These properties are independent of the choice of the admissible chart and the trivialization $E|_{U} \cong U \times \BC^{r}$.
It is clear that $(E,h)$ has a logarithmic divergence of order at most $k$ along $D$ if and only if its dual $(E^{*},h^{*}={}^{t}h^{-1})$ has a logarithmic decay of order at most $-k$ along $D$.

We also introduce a slightly weaker but more tractable notion, which is characterized by the divergence rate of the pointwise length of each section.

\begin{dfn}\label{LD}
Let $(X,D)$ and $(E,h)$ be as in \ref{adm}.
A singular Hermitian metric $h$ on $E$ is said to {\it have section-wise logarithmic divergence along $D$} if, for every $x \in D$, any admissible chart $(U,t=(t_{1},\ldots,t_{n}))$ centered at $x$, and any section $s \in H^{0}(U,E)$, there exist constants $c>0$ and $\delta \in (0,1)$ and a non-negative integer $k$ such that
$$ 
	|s|_{h}^{2}(t) \le c \prod_{j=1}^{n_{x}} (- \log|t_{j}|)^{k} 
$$ 
holds for any $t \in \{|t| < \delta\} \sm D$.
Here, the constants $c, \delta$, and $k$ may depend on $s$ (and on $x$).
\qedd
\end{dfn}

\begin{rem}
(1) If $h$ is Griffiths semi-negative, then $h$ is locally bounded from above.
That is, for each $x \in X$, there exist constants $c_{x}>0$ and $\delta_{x}>0$ such that $\lam_{r}(t) \le c_{x}$ for any $t \in \{|t|<\delta_{x}\} \sm D$.
In fact, the inequality $\lam_r(t) \le c_x$ holds even on $D$, though $\lam_1(t)$ may vanish on $D$.

(2) If $h$ is Griffiths semi-positive, then $h$ is locally bounded away from zero.
That is, for each $x \in X$, there exist constants $c_{x}>0$ and $\delta_{x}>0$ such that $\lam_{1}(t) \ge c_{x}$ for any $t \in \{|t|<\delta_{x}\} \sm D$.
In fact, the inequality $\lam_1(t) \ge c_x$ holds even on $D$, though $\lam_r(t)$ may become infinite on $D$.
\qedd
\end{rem}

We present a few related basic properties.
For the remainder of this section, we retain the notation $(X,D)$ and $(E,h)$ from \ref{adm}.

It is clear that condition \ref{LDD}\,(1) is stronger than \ref{LD}.
The converse also holds on every relatively compact open subset.

\begin{lem}\label{vsLD}
Suppose that a singular Hermitian metric $h$ satisfies \ref{LD}.
Then $h$ satisfies \ref{LDD}\,(1) on every relatively compact open subset.
\end{lem}

\begin{proof}
We choose a point $x \in D$ and an admissible chart $(U,t)$ centered at $x$.
We assume that $U = \{|t| < 1\} \subset \BC^{n}$ with coordinates $t=(t_{1},\ldots, t_{n})$, $x=0$, $D \cap U = \{ \prod_{j=1}^{n_{0}} t_{j}=0 \}$, and $E|_{U} \cong U \times \BC^{r}$.
We choose a local frame $e_{1},\ldots, e_{r} \in H^{0}(U, E)$ such that $E|_{U} = \bigoplus_{i=1}^{r} \CO_{U} e_{i}$.
As in \ref{adm}, we view $h$ as a matrix-valued function $h = ( h_{i\ol j} )$ on $U$.
By assumption \ref{LD}, after shrinking the neighborhood so that $-\log|t_{j}| \ge 1$ for all $j=1,\ldots,n_{0}$, we may choose constants $c_{0}>0$, $\delta_{0} \in (0,e^{-1})$, and $k_{0} \ge 0$, independent of $i$, such that for any $t \in \{|t| < \delta_{0} \} \sm D$,
$$
	|e_{i}|_{h}^{2}(t) \le c_{0} \prod_{j=1}^{n_{0}} (- \log|t_{j}|)^{k_{0}}
$$
holds for all $i=1, \ldots, r$.

Now, let $s \in H^{0}(U, E)$ be an arbitrary local section, and write $s = \sum_{i=1}^{r} s_{i}e_{i}$ with $s_{i} \in H^{0}(U,\CO_{U})$.
We define the standard Euclidean norm of $s$ by $|s|_{\mathrm{eu}} = (\sum_{i=1}^{r} |s_{i}|^{2})^{1/2}$.
Then, applying the Cauchy--Schwarz inequality, we obtain for each $t \in \{|t| < \delta_{0} \} \sm D$:
$$
	|s|_{h}^{2}(t) \le r \sum_{i=1}^{r} |s_{i}(t)|^{2} |e_{i}|_{h}^{2}(t)
	\le r c_{0} |s|_{\mathrm{eu}}^{2}(t) \prod_{j=1}^{n_{0}} (- \log|t_{j}|)^{k_{0}}.
$$
For any $t \in \{|t| < \delta_{0} \} \sm D$, we have
$$
	\lam_{r}(t) = \sup \Big\{ \sum_{i,j=1}^{r} h_{i \ol j}(t) v_{i}\ol{v_{j}} \mid
				v \in \BC^{r}, \ |v|_{\mathrm{eu}}=1 \Big\}.
$$
For each vector $v=(v_{1},\ldots,v_{r}) \in \BC^{r}$ with $|v|_{\mathrm{eu}}=1$, we associate a constant section $s_{v}= \sum_{i=1}^{r} v_{i}e_{i} \in H^{0}(U,E)$.
Since $|s_{v}|_{\mathrm{eu}} = 1$, it follows from the previous inequality that
$$
	\lam_{r}(t) = \sup \Big\{ |s_{v}|_{h}^{2}(t) \mid v \in \BC^{r}, \ |v|_{\mathrm{eu}}=1 \Big\}
		 \le  r c_{0} \prod_{j=1}^{n_{0}} (- \log|t_{j}|)^{k_{0}}.
$$
Hence, condition \ref{LDD}\,(1) holds on the neighborhood $U_{x} := \{|t| < \delta_{0}\}$ with the order at most $k_{0}$.

Now, let $W \subset X$ be a given relatively compact open subset. 
We can cover $\ol{W} \cap D$ by finitely many such neighborhoods $U_{x}$.
Away from $D$, the metric $h$ is locally bounded.
By setting $k_{W} = \max k_{x} \in \BZ_{\ge 0}$ over this finite collection, $h$ satisfies \ref{LDD}\,(1) on $W$ with the order at most $k_{W}$.
\end{proof}

Let us consider some functorial properties.

\begin{lem}\label{subq}
Let $(X,D)$ and $(E,h)$ be as in \ref{adm}.
Let $0 \to S \to E \to Q \to 0$ be a short exact sequence of holomorphic vector bundles on $X$.
The singular Hermitian metric $h$ on $E$ naturally induces singular Hermitian metrics $h_{S}$ on $S$ and $h_{Q}$ on $Q$, respectively.

(1) Suppose that $h$ has a logarithmic decay of order at most $-k$ along $D$.
Then $h_{S}$ also has a logarithmic decay of order at most $-k$ along $D$.

(2) Suppose that $h$ has a logarithmic divergence of order at most $k$ along $D$.
Then $h_{Q}$ also has a logarithmic divergence of order at most $k$ along $D$.
\end{lem}

\begin{proof}
(1) This is clear from the definition.
(2) This is reduced to (1) by taking the dual and considering the bundle injection $(Q^{*}, h_{Q}^{*}) \hookrightarrow (E^{*}, h^{*})$.
\end{proof}

The next lemma is also straightforward.

\begin{lem}\label{pullb}
Let $(X,D)$ and $(E,h)$ be as in \ref{adm}, and let $f: Y \to X$ be a smooth, surjective holomorphic map.
The preimage $D_{Y} = f^{-1}(D)$ is regarded as a normal crossing divisor on $Y$, while the pull-back bundle and metric are denoted by $(E_{Y}, h_{Y}) = (f^{*}E, f^{*}h)$.
If $h$ has a logarithmic divergence (resp.~logarithmic decay) of order at most $k$ (resp.~$-k$) along $D$, 
then $h_{Y}$ has a logarithmic divergence (resp.~logarithmic decay) of order at most $k$ (resp.~$-k$) along $D_{Y}$.
\end{lem}

A model case is the projective space bundle $\pi: \BP(E) \to X$ along with the natural quotient $\pi^*E \to \OE1$.
Hence, we obtain the following.
 
\begin{lem}\label{E2O1}
Let $(X,D)$ and $(E,h)$ be as in \ref{adm}. 
Suppose that $h$ is Griffiths semi-positive and satisfies \ref{LD}.
Then, the induced singular Hermitian metric $h_{\OE1}$ on $\OE1$ is Griffiths semi-positive and satisfies \ref{LD} along $\pi^{-1}(D) \subset \BP(E)$.
\end{lem}

The logarithmic divergence condition is sufficient for applying Demailly's criterion \ref{Dcor64}.

\begin{lem}
Suppose that $E$ is of rank 1 and that $h$ has semi-positive curvature and satisfies \ref{LD}.
Then, every local psh potential of $h$ has vanishing Lelong numbers everywhere. 
In particular, $E$ is metrically nef.
\end{lem}

\begin{proof}
In view of \ref{Dcor64}, the problem reduces to studying the singularities of $h$ on a local chart $(U,t)$ around a point $x \in D$.
We write $h = e^{-\vph}$, where $\vph$ is a psh function on $U$.
By virtue of \ref{vsLD}, our assumption implies that there exist constants $c_{x}>0$, $\delta_{x} \in (0,1)$, and $k_{x} \ge 0$ such that
$$
	e^{-\vph(t)} \le c_{x} \prod_{j=1}^{n_{x}} (- \log|t_{j}|)^{k_{x}}
$$
holds for any $t \in \{|t| < \delta_{x}\} \sm \{\prod_{j=1}^{n_{x}} t_{j}=0\}$, where $D \cap U = \{ \prod_{j=1}^{n_{x}} t_{j}=0\}$.
It follows that
$$
	\vph(t) \ge k_{x} \sum_{j=1}^{n_{x}} ( - \log (- \log|t_{j}|) ) - \log c_{x}. 
$$
This pointwise lower bound immediately implies that $\nu(\vph,x)=0$.
\end{proof}

In the setting of Theorem \ref{mt}, it therefore suffices to prove that the canonical $L^{2}$-metric on $E$ has section-wise logarithmic divergence along the boundary divisor $Y \sm U_{Y}$ on $Y$.

\begin{rem}
This is a side remark, and will not be referred to later.
Let $(X,D)$ and $(E, h)$ be as in \ref{adm}.
We can also consider various other classes of singularities for $h$, extending beyond logarithmic divergence and decay.
For example, when $X = \BC$ and $D = \{0\}$, one might consider the condition that 
$$
	|s|_{h}^{2}(t) \lesssim \frac1{|t|^{2a}} (-\log |t|)^{k}
$$
holds for $s \in H^{0}(U,E)$, where $a \in \BQ_{\ge 0}$ and $k \in \BZ_{\ge 0}$.

This type of singularity typically arises in a canonical $L^{2}$-metric when a fibration $f: X \to Y$ has singularities worse than semi-stable ones, such as non-reduced fibers.
In this setting where $Y = (\BC,t)$, the constant is given by $a = 1 - \mathrm{lct}(X, X_{0}) \in [0,1) \cap \BQ$, where $\mathrm{lct}(X, X_{0})$ denotes the log-canonical threshold of the pair $(X, X_{0})$ with the divisor $X_{0} = f^{*}(0)$.
\qedd
\end{rem}


\section{Logarithmic divergence of canonical $L^{2}$-metrics}

Let everything be as in \ref{mt}.
We will show that the direct image sheaf $E=f_{*}(K_{X/Y} \ot M)$ equipped with the canonical $L^{2}$-metric $g$ satisfies the properties in \ref{E2O1}.
Specifically, we will discuss (i) the definition of $g$, 
(ii) the Griffiths semi-positivity of $g$ (see \ref{Gsp}), and 
(iii) the property that $g$ has at worst logarithmic divergence (see \ref{LDn}).

We use the notation of \ref{def-wss} and set $Y_{0} = U_{Y}, X_{0}=U_{X}= f^{-1}(U_{Y})$, and $B = Y\sm Y_{0}$ the normal crossing divisor.
Set $E_{0}=E|_{Y_{0}}$ and $n= \dim X - \dim Y$.
Denote by $f_{0}:X_{0}\to Y_{0}$ the restriction of $f$, which is a smooth holomorphic map between complex manifolds.

We first define the canonical $L^{2}$-metric $g_{0}$ on $E_{0}$.
Let $V \subset Y_{0}$ be an open subset.
We set $X_{V}= f^{-1}(V)$.
We view a section
$$
	s \in H^{0}(V,E)= H^{0}(X_{V}, K_{X/Y}\ot M)
$$
as an $M$-valued relative holomorphic $n$-form on $X_{V}$.
We then consider
$$
	\w= h\, {\rm i}_n s\wed\ol s, \quad \text{and} \quad
	|s|_{g_{0}}^{2}(t) =\int_{X_{t}}\w
$$ 
for $t \in V$.
Here, ${\rm i}_n=\ai^{n^{2}}$ so that the pairing becomes semi-positive on each fiber.
By the theorem on differentiation under the integral sign, the function $|s|_{g_{0}}^{2} : V \to \BR_{\ge 0}=[0,+\infty)$ is $C^{\infty}$-smooth on $V$.
This defines a Hermitian metric $g_{0}$ on $E_{0}$ which is Nakano semi-positive by Berndtsson \cite{Be} (\cite{MT} when $M$ is a vector bundle).

\begin{lem}\label{Gsp}
The Hermitian metric $g_{0}$ on $E_{0}$ extends uniquely to a singular Hermitian metric $g$ on $E$ with Griffiths semi-positive curvature.
\end{lem}

\begin{proof}
We choose a point $0\in B \subset Y$.
We provide the unique extension of $g_{0}$ around 0 in $Y$.

We prepare a quasi-smooth model and a smooth model of $X$ to study $g_{0}$ and its extension $g$ to the whole space $Y$ (possibly restricting everything over any relatively compact open neighborhood of 0 in $Y$).
Let $a : X_{qs} \to X$ be a modification such that $f_{qs}=f\circ a :X_{qs}\to Y$ is a quasi-smooth weakly semi-stable model as in \ref{qs-wss}. 
In particular, $f_{qs}=f\circ a :X_{qs} \to Y$ is identical to $f:X\to Y$ over $Y_{0}$.
(When $X$ has at worst quotient singularities in \ref{mt}, we take $X_{qs}=X$ and $a : X_{qs}\to X$ to be the identity map.)
Put $M_{qs}=a^{*}M$ and $h_{qs}=a^{*}h$ which is Nakano semi-positive.
We further take a resolution of singularities $b:X_{sm} \to X_{qs}$ such that $X_{sm}$ is smooth, $b$ is projective, and the composition $f_{sm}=f\circ a \circ b :X_{sm}\to Y$ is identical to $f:X\to Y$ over $Y_{0}$.
Put $M_{sm}=(a\circ b)^{*}M$ and $h_{sm}=(a\circ b)^{*}h$ which is Nakano semi-positive.
We can take $a$ and $b$ so that $f_{sm}$ is locally K\"ahler (see \cite[6.2]{takegoshi}).

Let us observe that the direct image sheaf $E= f_{*}(K_{X/Y}\ot M)$ is independent of the models $f_{qs}$ and $f_{sm}$, namely
$$
	(f_{sm})_{*} (K_{X_{sm}/Y}\ot M_{sm}) 
	= (f_{qs})_{*} (K_{X_{qs}/Y}\ot M_{qs})
	=E.
$$
Since $X$ has at worst canonical singularities, we have $(a\circ b)_{*}K_{X_{sm}} = K_{X}$ and $a_{*}K_{X_{qs}} = K_{X}$.
Moreover, it follows that
$$
\begin{aligned}
	(f_{sm})_{*}(K_{X_{sm}/Y} \ot M_{sm})
	& = f_{*}(a \circ b)_{*}(K_{X_{sm}} \ot (a\circ b)^{*} (f^{*}K_{Y}^{-1} \ot M))
	\\
	& = f_{*} ((a\circ b)_{*}K_{X_{sm}} \ot f^{*}K_{Y}^{-1} \ot M)
	=  f_{*}(K_{X/Y} \ot M).
\end{aligned}
$$
Since $X_{sm}$ is smooth and $(M_{sm},h_{sm})$ is Nakano semi-positive, we can directly apply the known results to $f_{sm}: X_{sm} \to Y$ (see, e.g., \cite{BP}, \cite{PT}, \cite{HPS}, \cite{Cao}, \cite{ZZ}, and \cite{Ta24}).
The last three references treat the K\"ahler setting.
More precisely, the Ohsawa-Takegoshi type optimal $L^{2}$-extension theorem in the K\"ahler setting from \cite{Cao} (when $M$ is of rank 1) and \cite{ZZ} can be applied to the argument of \cite{HPS}, which assumes the projectivity of the family.
Consequently, the direct image $E_{sm} := (f_{sm})_{*}(K_{X_{sm}/Y} \ot M_{sm})$ is Griffiths semi-positive with respect to the canonical $L^{2}$-metric $g_{sm}$.
As a consequence of the theorems mentioned above, $g_{sm}$ is defined on $E_{sm}=E$ as a singular Hermitian metric globally, not just on $E_{0}$.

Since $f_{sm}$ and $f$, as well as $(M_{sm},h_{sm})$ and $(M,h)$, coincide over $Y_{0}$, we have $(E_{sm}, g_{sm})|_{Y_{0}} = (E_{0},g_{0})$.
Hence, $g_{0}$ extends to $g_{sm}$ on $E$.
The uniqueness of such extensions follows from the lemma below.
\end{proof}

The uniqueness of extensions is a more abstract property, which is easier to examine in the dual setting.

\begin{lem}
Let $F$ be a holomorphic vector bundle on a complex manifold $Z$.
Let $h_{1}$ and $h_{2}$ be singular Hermitian metrics on $F$ with Griffiths semi-negative curvature.
Suppose that $h_{1}=h_{2}$ on a non-empty Zariski open subset $Z_{0}$ of $Z$.
Then $h_{1}=h_{2}$ holds everywhere on $Z$.
\end{lem}

\begin{proof}
The assertion is local on $Z$.
We choose a small open subset $U \subset Z$ and an arbitrary local section $s \in H^{0}(U,F)$.
Since the curvature is Griffiths semi-negative, both $\log |s|_{h_{i}}$ and $|s|_{h_{i}}$ are psh on $U$ for each $i=1,2$.
By assumption, we have $|s|_{h_{1}}=|s|_{h_{2}}$ on $U \cap Z_{0}$.
Since $U \cap Z_{0}$ is a dense Zariski open subset of $U$, the psh functions $|s|_{h_{1}}$ and $|s|_{h_{2}}$ agree almost everywhere on $U$ and hence everywhere on $U$.
Since $s$ is arbitrary, this immediately leads to $h_{1} = h_{2}$ on $U$, completing the proof.
\end{proof}

We now have a singular Hermitian metric $g$ on $E$ with Griffiths semi-positive curvature.
To construct $g$ from $g_{0}$ and establish its Griffiths semi-positivity, we utilized the resolution $f_{sm}: X_{sm} \to Y$.
Although it would be possible to obtain $g$ by working directly with $f$, doing so would require a more explicit computation of fiber integrals.

We shall use the quasi-smooth model to prove that $g$ has at worst logarithmic divergence along $B=Y\sm Y_{0}$.
This is indispensable for estimating the singularities of $g$.
For this purpose, we may suppose that $X$ is quasi-smooth from the beginning.
By combining the discussions above, the following theorem completes the proof of \ref{mt}.

\begin{thm}\label{LDn}
Suppose in \ref{mt} that $X$ is quasi-smooth.
Then, the canonical $L^{2}$-metric $g$ on $E$ has a logarithmic divergence of order at most $n$ along $B$, where $n = \dim X - \dim Y$.
\end{thm}

\begin{proof}
Let $\dim Y = m$ and $\dim X = m+n$.
We choose a point $y_{0} \in B \subset Y$ and an admissible chart $(V,t) \subset Y$ centered at $y_{0}$ such that $B \cap V = \{\prod_{j=1}^{m_{y_{0}}} t_{j} = 0\}$ as in \ref{adm}.
Then, the proof of \ref{LDn} reduces to showing the following claim.

\smallskip

\noindent
{\bf Claim.}
{\it Let $s \in H^{0}(V,E)$.
Then, there exist constants $c>0$ and $\delta \in (0,1)$ such that, for any $t \in \{|t| < \delta\} \sm B$,
$$ 
	|s|_{g}^{2}(t) \le c \prod_{i=1}^{m_{y_{0}}} (-\log|t_{i}|)^{n}
$$ 
holds.}

\smallskip

This is more or less known to experts; however, we provide a proof for the sake of completeness.
As we will see below, the divergence of $|s|_g^2$ comes only from the singularities of $f$.
In fact, when we write $s$ as a relative differential form, it has singularities only where $f$ is not smooth.

We set $X_{V}= f^{-1}(V)$.
We view 
$
s \in H^{0}(V,E)= H^{0}(X_{V}, K_{X/Y}\ot M)
$
as an $M$-valued relative holomorphic $n$-form on $X_{V}$.
We then consider 
$
	\w=h\, {\rm i}_n s\wed\ol s
$
and 
$$
	|s|_{g}^{2}(t) =\int_{X_{t}}\w
$$ 
for $t \in V\sm B$.
The function $|s|_{g}^{2} : V \sm B \to \BR_{\ge 0}=[0,+\infty)$ is $C^{\infty}$-smooth on $V\sm B$.
We now analyze the behavior of $|s|_{g}^{2}(t)$ as $t \to 0$.

\medskip
\noindent
{\it Part I: Semi-stable case.} 
We first assume that $X$ is smooth.

{\it Step 1: Localization.}
We choose a locally finite, sufficiently fine open covering $\{U_{\lam}\}$ of $X$, and choose a $C^{\infty}$-smooth partition of unity $\{\rho_{\lam}\}$ subordinate to $\{U_{\lam}\}$.
Then, for every $t \in V\sm B$, we have
$$
	|s|_{g}^{2}(t) = \sum_{\lam} \int_{X_{t}\cap U_{\lam}} \rho_{\lam} \w,
$$
where the summation is actually a finite sum because each fiber $X_{t}$ is compact.
After shrinking $V$, properness of $f$ and local finiteness of the covering allow us to work with finitely many charts uniformly for $t$ near $y_{0}$.
Hence, the study of the singularities or the asymptotic behavior of $|s|_{g}^{2}(t)$ around $y_{0} \in B$ reduces to the local computation of $\int_{X_{t}\cap U_{\lam}} \rho_{\lam} \w$ on each $U_{\lam}$.
We denote the central fiber by $X_{y_{0}} = f^{-1}(y_0)$.

\smallskip
{\it Step 2: Local coordinates.}
For each $x_{0} \in X_{y_{0}}$, there exist local coordinates on $U \subset \BC^{n+m}$ denoted by $x=(x_{1},\ldots,x_{n+m})$ and on $V \subset \BC^{m}$ denoted by $t=(t_{1},\ldots,t_{m})$, such that the map $f|_{U}: U \to V$ is given by $t=f(x)$, $t_{j}=f_{j}(x)$ for $j=1,\ldots,m$, with
$$
	t_{1}=\prod_{j=1}^{\ell_{1}} x_{j}, \quad
	t_{2}=\prod_{j=\ell_{1}+1}^{\ell_{2}} x_{j}, \quad \ldots, \quad 
	t_{m}=\prod_{j=\ell_{m-1}+1}^{\ell_{m}} x_{j} 
$$
for some integers $0=\ell_{0}<\ell_{1}<\ldots < \ell_{m}\le n+m$.
(The coordinate systems $x$ and $t$ depend on $x_{0}$, and the initial $(V,t)$ may need to be replaced.
However, this does not affect the desired estimate.)
We can write $\ell_{i}=\ell_{i-1}+n_{i}+1$ with $n_{i} \in \BZ_{\ge 0}$ and $\sum_{i=1}^{m+1} n_{i}=n$, where $n_{m+1} := n+m-\ell_{m}$.
We denote $y_{j}=x_{\ell_{j}}$ for $j=1,\ldots,m$, and set
$$ 
\begin{aligned}
&\bfx_{j}=(x_{\ell_{j-1}+1},\ldots,x_{\ell_{j-1}+n_{j}}), 
\quad
&\bfx_{m+1}=(x_{\ell_{m}+1},\ldots,x_{n+m}).
\end{aligned}
$$ 
We use the following conventions.
For $j=1,\ldots,m$, set 
$$
	\xi_j :=\prod_{i=\ell_{j-1}+1}^{\ell_{j-1}+n_j}x_i, 
	\quad
	d\bfx_j := \bigwedge_{i=\ell_{j-1}+1}^{\ell_{j-1}+n_j} dx_i,
	\quad
	d\ol \bfx_j := \bigwedge_{i=\ell_{j-1}+1}^{\ell_{j-1}+n_j} d \ol x_i.
$$
Then $t_j=\xi_j y_j$ for $j=1,\ldots,m$.
We set $d\bfx_{m+1}$ and $d\ol\bfx_{m+1}$ in the same way.
Empty products and empty exterior products are understood to be $1$.

We can still assume $B \cap V = \{ \prod_{j=1}^{m_{y_{0}}} t_{j} = 0 \} $ for some $1 \le m_{y_{0}} \le m$.
For notational simplicity, we assume $m_{y_{0}} = m$; the other cases where $m_{y_{0}} \le m-1$ can be easily deduced from the case $m_{y_{0}} = m$.

For $t=(t_{1},\ldots,t_{m}) \in V$, we set
$$
	X_{t_{j}} = \{(\bfx_{j}, y_{j}) \in \BC^{n_{j}+1} \mid t_{j}= \xi_j y_{j} \}
$$ 
for $j=1,\ldots,m$, and set $X_{m+1} = \{\bfx_{m+1} \in \BC^{n+m-\ell_{m}} \} = \BC^{n_{m+1}}$.
Then, the fiber of $f|_{U}: U \to V$ is given by the product
$$
	X_{t} = X_{t_{1}} \times \ldots \times X_{t_{m}} \times X_{m+1}
$$
in $\BC^{n+m} = \BC^{n_{1}+1} \times \ldots \times  \BC^{n_{m}+1} \times \BC^{n_{m+1}}$.
(Although the original space $X$ is not identical to the entire $\BC^{n+m}$ globally, this does not affect the estimation of the integrals over compact subsets.)

\smallskip

{\it Step 3: Integrand.}
Let $dx = \bigwedge_{i=1}^{n+m} dx_{i}$ and $dt = \bigwedge_{j=1}^{m} dt_{j}$ be local frames of $K_{X}|_{U}$ and $K_{Y}|_{V}$, respectively.
Using the relation $\frac{dt_{j}}{t_{j}} = \sum_{i=\ell_{j-1}+1}^{\ell_{j-1}+n_{j}} \frac{dx_i }{ x_{i} } + \frac{dy_{j}}{y_{j}}$, we obtain
$$
	\bigwedge_{i=\ell_{j-1}+1}^{\ell_{j}} dx_{i}
	=  \frac{y_{j}}{t_{j}} \Big( \bigwedge_{i=\ell_{j-1}+1}^{\ell_{j-1}+n_{j}} dx_{i} \Big) \wed dt_{j}.
$$
In other words, this can be written as $d\bfx_{j} \wed dy_{j} = \frac{y_{j}}{t_{j}} d\bfx_{j} \wed dt_{j} = \frac{ d\bfx_{j} }{ \xi_j } \wed dt_{j}$.
Then, a local frame $\frac{dx}{dt}$ of $K_{X/Y}$ on $U$ can be expressed as
$$
	\frac{dx}{dt} 
	= \bigwedge_{j=1}^{m} 
	\Big(\frac{y_{j}}{t_{j}}  d\bfx_{j}\Big) \wed d\bfx_{m+1}
	= \bigwedge_{j=1}^{m}
	 \frac{ d\bfx_{j} }{ \xi_j} \wed d\bfx_{m+1},
$$
which satisfies $\frac{dx}{dt} \wed f^{*}dt  = \pm dx$.
Regarding $s$ as an element of $H^{0}(X_{V},K_{X/Y} \ot M) \cong \mathrm{Hom}_{\CO_{X_{V}}}(f^{*}K_{Y}, K_{X} \ot M)$, we have $s(f^{*}dt)|_{U} = s_{U} dx = \pm s_U \frac{dx}{dt} \wed f^{*}dt$ for some vector-valued holomorphic function $s_{U} \in H^{0}(U,\CO_{X}^{\oplus r})$ under a local trivialization of $M$.
Hence, it follows that
$$
	s|_{X_{t}} = \pm\, s_U \frac{dx}{dt} \Big|_{X_{t}}
$$ 
for $t \in V \sm B$.
We can represent $h=(h_{j\ol k})$ as a continuous (in fact, smooth) positive-definite matrix-valued function on $U$, and $s_U={}^t(s_{U1},\ldots,s_{Ur})$ as a holomorphic column-vector-valued function on $U$.
Then 
$$
	h\, {\rm i}_n s \wed \ol s|_{U \cap X_t}
	= \Big(\sum_{j,k} h_{j\ol k}s_{Uj} \ol{s_{Uk}} \Big) \
		{\rm i}_n \frac{dx}{dt} \Big|_{X_{t}} \wed \ol{\frac{dx}{dt} }\Big|_{X_{t}}.
$$
We write $h|s_U|^2 := \sum_{j,k} h_{j\ol k}s_{Uj} \ol{s_{Uk}}$.
Consequently, the integrand $\rho \omega$ restricted to $U \cap X_{t}$ is given by
$$
	\rho h |s_{U}|^{2} \bigwedge_{j=1}^{m+1} W_{j},
$$
where
$$
	 W_{j} = \frac{{\rm i}_{n_{j}} d\bfx_{j} \wed d\ol\bfx_{j}  }{ |\xi_{j}|^{2} } \quad (j=1,\ldots,m) 
	 \quad \text{and} \quad
	 W_{m+1} = {\rm i}_{n_{m+1}} d\bfx_{m+1} \wed d\ol\bfx_{m+1}.  
$$

\smallskip

{\it Step 4: Estimate.}
We will estimate the fiber integral by applying Fubini's theorem.
Let $R: [0,+\infty) \to \BR_{\ge 0}$ be a $C^{\infty}$-smooth cut-off function such that $R \equiv 1$ on $[0,c_{0}]$ and $R \equiv 0$ on $[c_{0}+1,+\infty)$ for a given constant $c_{0}>0$.
We define $\wtil R(x) := \prod_{i=1}^{n+m} R(|x_{i}|) \in C^{\infty}_{\cpt}(\BC^{n+m},\BR_{\ge 0})$ and set $\wtil U = \prod_{i=1}^{n+m}\{ |x_{i}| < c_{0}+1\}$.
Here, the constant $c_{0}>0$ can be chosen sufficiently large so that $\wtil R \equiv 1$ on $\Supp \rho$ for the pre-assigned cut-off function $\rho$ on $U$ (such as $\rho_{\lam}$ in Step 1).
Then, there exists a constant $A>0$ such that $0 \le \rho h |s_{U}|^{2} \le A \cdot \wtil R(x)$ holds on $U$.
We can express the dominating volume form as
$$
	\wtil R(x) \bigwedge_{j=1}^{m+1} W_{j}
	= \bigwedge_{j=1}^{m+1} \wtil W_{j},
$$
where
$$
	\wtil W_{j} = R(|y_{j}|) \bigwedge_{i=\ell_{j-1}+1}^{\ell_{j-1}+n_j} R(|x_{i}|) \frac{  \ai dx_{i} \wed d\ol x_{i}  } { | x_{i} |^{2}  } \quad (j=1,\ldots,m),
$$
and $\wtil W_{m+1} = \bigwedge_{i=\ell_m+1}^{n+m} R(|x_{i}|) \ai dx_{i} \wed d\ol x_{i}$.
We then obtain
$$
	0 \le \int_{X_{t}\cap U}\rho \w
	\le A \Big( \prod_{j=1}^{m} I_{j} \Big) \times I_{m+1},
$$
where
$$ 
\begin{aligned}
	I_{j} = I_{j}(t_{j}) 
	&= \int_{ \BC^{n_j} \sm \{\xi_{j}=  0\}} 
 	R\Big( \Big| \frac{t_j}{ \xi_{j} }  \Big| \Big) 
 	\bigwedge_{i=\ell_{j-1}+1}^{\ell_{j-1}+n_j} R(|x_{i}|) 
 	\frac{  \ai dx_{i} \wed d\ol x_{i}  } { | x_{i} |^{2}  },
\\
I_{m+1} & =
 \int_{ \BC^{n_{m+1}} }
 \bigwedge_{i=\ell_m+1}^{n+m} R(|x_{i}|) \ai dx_{i} \wed d\ol x_{i}.  
\end{aligned}
$$
The last integral $I_{m+1}$ over $\BC^{n_{m+1}}$ is a finite positive number.

Each integral $I_{j}(t_{j})$ for $j = 1,\ldots,m$ can be computed by introducing polar coordinates in each variable.
We illustrate this calculation explicitly for $I_{1}$.
Let $r=(r_{1},\ldots,r_{n_{1}}) \in (\BR_{\ge 0})^{n_{1}}$ and $\theta= (\theta_{1},\ldots,\theta_{n_{1}}) \in [0,2\pi]^{n_{1}}$ with $x_{i} = r_{i}e^{\ai \theta_{i}}$ for $i=1,\ldots,n_{1}$.
Then, we have
$$
	 I_{1}(t_{1}) = \int_{r,\theta} 
	 R\Big( \frac{|t_{1}|}{  \prod_{i=1}^{n_{1}} r_{i} }  \Big) 
	 \prod_{i=1}^{n_{1}} R( r_{i} ) 
	 \frac{  2 dr_{i} d\theta_{i}  } {  r_{i}  }.
$$
By setting $s = |t_{1}| > 0$, the integral can be rewritten as
$$
	 I_{1}(t_{1}) = F_{1}(s) := 2^{n_{1}} (2\pi)^{n_{1}} 
	 \int_{ ( \BR_{\ge 0} )^{n_{1}} } 
	  R \Big( \frac{ s }{  \prod_{i=1}^{n_{1}} r_{i} }  \Big) 
	 \prod_{i=1}^{n_{1}} R( r_{i} ) r_{i}^{-1}\, dr_{i}.
$$
According to \cite[Prop.\,1\,(p.\,131)]{Ba}, there exist constants $p_{k} \in \BR$ and $\ep_{1} \in (0,1)$ such that
$$
	F_{1}(s) = \sum_{k=0}^{n_{1}} p_{k} (-\log s)^{k}
$$
for all $s \in (0,\ep_{1})$.
In particular, there exist constants $c_{1}>0$ and $\del_{1} \in (0,\ep_{1})$ such that the bound $0 \le I_{1}(t_{1}) \le c_{1}(-\log |t_{1}|)^{n_{1}}$ holds for any $t_{1} \in \{|t_{1}|<\del_{1} \} \sm \{0\}$.

In summary, we have shown that there exists a constant $A'>0$ such that 
$$
	\int_{X_{t}\cap U} \rho \w
	\le A' \prod_{j=1}^m (-\log |t_j|)^{n_{j}}
$$
holds for every $t \in V \sm B$, possibly after shrinking $V$.
Here, the exponents satisfy $n_{j} \le n$ for all $j$.

\medskip
\noindent
{\it Part II: General case.}
Suppose that $f: X \to Y$ is weakly semi-stable and that $X$ is quasi-smooth.

The problem can be reduced to the semi-stable case as follows.

{\it Step 1: Localization.}
For every $x \in X$, there exists an open neighborhood $U \subset X$ of $x$ such that $U = \tU/G$, where $\tU \subset \BC^{m+n}$ is an open neighborhood of the origin and $G$ is a finite abelian group depending on $x$.
Furthermore, the composition $f \circ \pi$ is semi-stable, which means it is given by monomials whose exponents are all equal to one, as described in Part I, Step 2, where $\pi: \tU \to U$ denotes the quotient map.
We cover $X$ by such open subsets $\{ U_{\lam} \}$, each equipped with a quotient map $\pi_{\lam}: \tU_{\lam} \to U_{\lam}$ and a representation $U_{\lam} = \tU_{\lam}/G_{\lam}$.
Let $\{\rho_{\lam}\}$ be a partition of unity subordinate to $\{ U_{\lam} \}$.
For any $t \in V \sm B$, the integral to be estimated is expressed as
$$ 
\begin{aligned}
	|s|_{g}^{2}(t)
	& =\sum_{\lam} \int_{X_{t} \cap U_{\lam} } \rho_{\lam} \w, \\ 
	\int_{X_{t} \cap U_{\lam} } \rho_{\lam} \w 
	&  = \frac1{|G_{\lam}|} \int_{ (f\circ \pi_{\lam})^{-1}(t) } \pi_{\lam}^{*} (\rho_{\lam}\w).
\end{aligned}
$$
We have already estimated the last integral
$$
	\int_{ (f\circ \pi_{\lam})^{-1}(t) } \pi_{\lam}^{*} (\rho_{\lam}\w)
$$
in Part I. 

Repeating Steps 2--4 of Part I yields the desired estimate.

\medskip
This establishes the Claim, and hence completes the proof of \ref{LDn}.
\end{proof}


\section{Around Nakano semi-positivity}

We discuss metric positivity of vector bundles on possibly singular complex spaces in order to make our statements precise.
Readers who are familiar with these notions may skip this section.
We consider two approaches; see \ref{Nsp} and \ref{pushw}.

\smallskip

Let $X$ be a reduced and irreducible complex space. 
Let $E$ be a holomorphic vector bundle on $X$ equipped with a Hermitian metric $h$. 
Note that the metric $h$ is defined on the whole of $X$, rather than just on the regular part $X_{reg}$. 
Here, we adopt a standard notion of Hermitian metrics on complex spaces. 
Let $\{U_{i}\}_{i}$ be an open covering of $X$ such that each $U_{i}$ admits a local holomorphic embedding into a domain in some $\BC^{N}$ and a local trivialization of $E$. 
A Hermitian metric $h$ on $E$ is given by a collection of smooth matrix-valued functions $\{h_i\}_i$ satisfying the usual transition rule, where smoothness is understood in the sense described in \cite[\S 1]{fujiki}, for example.
In what follows, we will focus on the Nakano semi-positive case; however, the same arguments work for the Griffiths semi-positive case as well.

\begin{dfn}\label{Nsp}
(1)
We say that a Hermitian metric $h$ on $E$ is {\it Nakano semi-positive at $x \in X$} if there exists an open neighborhood $U$ of $x$ and a surjective holomorphic map $f : V \to U$ from a complex manifold $V$ such that the pull-back bundle $f^{*}E$ equipped with the pull-back metric $f^{*}h$ is Nakano semi-positive (in the usual sense). 

(2) 
We say that $h$ is {\it Nakano semi-positive} if it is Nakano semi-positive at every point $x \in X$.
Equivalently, $h$ is Nakano semi-positive on $X_{reg}$ in the usual sense (see \ref{smooth on $X_{reg}$}).

(3) We say that $E$ is {\it Nakano semi-positive} if it admits a Nakano semi-positive Hermitian metric.
\qedd
\end{dfn}

For example, in \ref{Nsp}, if $x \in X$ is a quotient singularity, we may take $f:V \to U$ to be the quotient map $\pi: \BC^{n} \to \BC^{n}/G$ by a finite group $G$ as a local model.
This definition \ref{Nsp} does not depend on the choice of the surjective holomorphic map $f:V \to U$ from a complex manifold $V$.
For example, we can take $f:V\to U$ as a resolution of singularities of $U$.
The independence of this definition from the choice of the surjective map follows from \ref{any-resol}--\ref{any-pb}.

\begin{lem}\label{any-resol}
Suppose that $X$ is a complex manifold (i.e., $X$ is smooth).
Let $f:Y \to X$ be a surjective holomorphic map from a complex manifold $Y$.
Then, a Hermitian metric $h$ on $E$ is Nakano semi-positive (in the usual sense) if and only if $f^{*}E$ equipped with the metric $f^{*}h$ is Nakano semi-positive (in the usual sense).
\end{lem}

\begin{proof}
Suppose that $f^{*}h$ is Nakano semi-positive.
Let $W \subset X$ be the set of points $x \in X$ for which there exists a point $y \in f^{-1}(x)$ such that the differential $df_{y} : T_{Y,y} \to T_{X,x}$ is surjective.
By assumption, for any $y \in Y$, the curvature tensor $(f^{*}\Theta_{h})(y)$ is Nakano semi-positive on $T_{Y,y} \otimes (E_{f(y)},h_{f(y)})$.
Thus, at each point $x \in W$, the tensor $\Theta_{h}(x)$ is Nakano semi-positive on $T_{X,x} \otimes (E_{x},h_{x})$.
Since the set $W$ is everywhere dense in $X$ and $\Theta_{h}$ is continuous, we conclude that $\Theta_{h}$ is Nakano semi-positive on the whole of $X$.

The converse direction is clear.
\end{proof}

The following alternative characterization is also useful.

\begin{lem}\label{smooth on $X_{reg}$}
The Hermitian metric $h$ on $E$ is Nakano semi-positive in the sense of \ref{Nsp} if and only if $h$ is Nakano semi-positive on the regular locus $X_{reg}$ in the usual sense.
\end{lem}

\begin{proof}
Nakano semi-positivity in the sense of \ref{Nsp} clearly implies Nakano semi-positivity on $X_{reg}$.
Let us show the converse.

For any $x \in X$, we choose a relatively compact open neighborhood $U \subset X$ of $x$ and a resolution of singularities $f:V \to U$ such that $f$ is biholomorphic over the regular locus $U_{reg}$.
Suppose that $h$ is Nakano semi-positive on $U_{reg}$.
Then, $f^{*}h$ is Nakano semi-positive on $f^{-1}(U_{reg})$.
Since $f^{*}h$ is smooth and $f^{-1}(U_{reg})$ is everywhere dense in $V$, it follows that $f^{*}h$ is Nakano semi-positive on the whole of $V$.
Hence, $h$ is Nakano semi-positive at $x$ in the sense of \ref{Nsp}.
\end{proof}

The property \ref{Nsp} is preserved by arbitrary pull-backs.

\begin{lem}\label{any-pb}
Suppose that a Hermitian metric $h$ on $E$ is Nakano semi-positive in the sense of \ref{Nsp}.
Let $g : Z \to X$ be a holomorphic map from a reduced and irreducible complex space $Z$.
Then, $g^{*}E$ equipped with the metric $g^{*}h$ is Nakano semi-positive in the sense of \ref{Nsp}.
\end{lem}

\begin{proof}
Fix a point $z \in Z$.
Our assertion is local around $z$.
We set $x = g(z) \in X$.
We can take an open neighborhood $U$ of $x$ and a surjective holomorphic map $f : V \to U$ from a complex manifold $V$ as in \ref{Nsp}.
Let $Z_{U} := g^{-1}(U)$ be an open neighborhood of $z$.

By considering the fiber product $Z_{U} \times_{U} V$ and taking its resolution of singularities, we can find a complex manifold $W$ along with a surjective holomorphic map $f' : W \to Z_{U}$ and a holomorphic map $g' : W \to V$ such that $g \circ f' = f \circ g'$ on $W$.
We then consider the pull-back of $E$ and $h$ to $W$, denoted by $E_{W}$ and $h_{W}$, respectively.
Since $h_{V} := f^{*}h$ is Nakano semi-positive on $V$, its pull-back $h_{W} = {g'}^{*}h_{V}$ is also Nakano semi-positive on $W$.
This implies that ${f'}^{*}(g^{*}h) = h_{W}$ is Nakano semi-positive.
Therefore, by definition, $g^{*}h$ is Nakano semi-positive at $z$.
\end{proof}

Next, with the aim of generalizing \ref{mt}, we introduce a slightly weaker version of Nakano semi-positivity, which is defined as Nakano semi-positivity after a certain pull-back of the vector bundle (rather than its Hermitian metric).

\begin{prop}\label{pushw}
Let $X$ be a reduced and irreducible complex space, and $E$ a holomorphic vector bundle on $X$.
Let $Y$ be a complex manifold, and $f:Y \to X$ a proper surjective holomorphic map with connected fibers.
Suppose that the pull-back bundle $f^{*}E$ admits a Hermitian metric $g$ with Nakano (resp.~ Griffiths) semi-positive curvature.
Then, $g$ is constant on each fiber of $f$.
The push-forward $h:=f_{*}g$ defines a continuous and everywhere positive definite Hermitian metric on $E$.
Moreover, $h=f_{*}g$ is smooth and Nakano (resp.~ Griffiths) semi-positive on the open subset of $X$ over which $f$ is smooth.
\end{prop}

A $C^{0}$-Hermitian metric can be viewed as a singular Hermitian metric in the sense that it lacks smoothness. 
When $X$ is normal, the push-forward $f_{*}g$ is in fact a Griffiths semi-positive singular Hermitian metric on $E$. 
(However, the notions of Griffiths semi-positivity and semi-negativity for singular Hermitian metrics are not well-developed on non-normal complex spaces.)
Moreover, $f_{*}g$ might be expected to be a Nakano semi-positive singular Hermitian metric on $E$; however, it remains unclear whether the notion of singular Nakano semi-positivity has been firmly established.
This depends on the definition of singular Nakano semi-positivity being used.

The property that $f_{*}g$ is $C^{0}$ is sufficient to give an estimate of fiber integrals from above in variants of \ref{mt} (see \ref{mt-variant}). 
A model case is when $X$ is normal and $f : Y \to X$ is a bimeromorphic holomorphic map from a complex manifold $Y$.

\begin{proof}[Proof of \ref{pushw}]
Since the assertion is local on $X$, we may assume that $E = X \times \BC^{r}$.
Let $0 \in X$ be a reference point.

We choose Hermitian metrics $\w_{X}$ on $X$ and $\w_{Y}$ on $Y$, respectively.
Let $d_{X}$ and $d_{Y}$ denote the induced distance functions, respectively.
For every $c > 0$, we define $U_{c} := \{ x \in X \mid d_{X}(0,x) < c \}$.
By choosing an appropriate metric $\w_{X}$ (e.g., realizing $X$ as a closed complex subspace of the unit ball $B$ of some $\BC^{N}$ and taking $\w_{X} = \w_{B}|_{X}$, where $\w_{B}$ is the Poincar\'e metric on $B$), we may assume that $U_{1}$ is relatively compact in $X$.
We set $Y_{1} := f^{-1}(U_{1})$, which is relatively compact in $Y$.
Let $M:=M_r(\BC)\cong\BC^{r^{2}}$ be the space of $r\times r$ complex matrices.

(1)
Lemma \ref{trivial} below shows that $g$ is constant on each fiber $Y_{x}$ (endowed with the induced reduced structure); namely, $g|_{Y_{x}} = Y_{x} \times h_{x}$ for some positive definite Hermitian matrix $h_{x} \in M$.
(For this step, singular Griffiths semi-positivity suffices as long as $g|_{Y_{x}}$ is well-defined.)

(2)
We shall prove that the map $h : X \to M$ given by $x \mapsto h_{x}$ is continuous at $x = 0$.

Let $d_{M}$ be the Euclidean distance function on $M=\BC^{r^{2}}$.
Let $\varepsilon > 0$.
Since $g$ is continuous and $\overline{Y_{1}}$ is compact, $g$ is uniformly continuous on $\overline{Y_{1}}$. Thus, there exists a constant $c_{g} > 0$ such that for any $y_{1}, y_{2} \in Y_{1}$ satisfying $d_{Y}(y_{1}, y_{2}) < c_{g}$, we have $d_{M}(g(y_{1}), g(y_{2})) < \varepsilon$.
We consider the following neighborhood of $Y_{0}$:
$$
V_{c_{g}} := \{ y \in Y \mid d_{Y}(y, Y_{0}) < c_{g} \},
$$
where we define $d_{Y}(y, Y_{0}) := \inf_{p \in Y_{0}} d_{Y}(y, p)$ as usual.
As will be shown in (3) below, $f(V_{c_{g}})$ contains an open neighborhood of $0 \in X$.
Let $x \in U_{1} \cap f(V_{c_{g}})$, and choose $\delta \in (0,1)$ such that $U_{\delta} \subset f(V_{c_{g}})$.
Then, we can take $y \in V_{c_{g}}$ such that $f(y) = x$.
By definition of $V_{c_{g}} $, we can also find a point $y_{0} \in Y_{0}$ satisfying $d_{Y}(y, y_{0}) < c_{g}$.
This implies that $d_{M}(g(y), g(y_{0})) < \varepsilon$, which means $d_{M}(h_{x}, h_{0}) < \varepsilon$.
Consequently, for any $x \in U_{\delta}$, we obtain $d_{M}(h_{x}, h_{0}) < \varepsilon$.

(3)
For every $c > 0$, we set $V_{c} := \{ y \in Y \mid d_{Y}(y, Y_{0}) < c \}$, which is an open neighborhood of $Y_{0}$.
We shall show that for every $c > 0$, $f(V_{c})$ contains a neighborhood of $0 \in X$.

Suppose that $f(V_{c})$ does not contain any neighborhood of $0 \in X$.
Then, there exists a sequence of points $\{x_{k}\}_{k \in \mathbb{N}}$ in $U_{1}$ such that $d_{X}(x_{k}, 0) \to 0$ as $k \to \infty$, but $x_{k} \notin f(V_{c})$ for all $k$.

For each $k$, we choose a point $y_{k} \in Y_{x_{k}} \subset Y_{1}$.
Since $\overline{Y_{1}}$ is compact, by passing to a subsequence if necessary, we may assume that the sequence $\{y_{k}\}$ converges to some point $q \in \overline{Y_{1}}$.
The limit $q$ must lie in $Y_{0}$ because $f(q) = f(\lim y_{k}) = \lim f(y_{k}) = \lim x_{k} = 0$.
Hence, we have $y_{k} \in V_{c}$ for all sufficiently large $k$.
This implies that $x_{k} = f(y_{k}) \in f(V_{c})$, which yields a contradiction.

(4)
Suppose that $f$ is smooth at $y \in Y$.
Set $x = f(y)$ and choose a sufficiently small open neighborhood $U$ of $x$.
Then, we can find a local section $s : U \to f^{-1}(U)$ of $f$ such that $s(x) = y$.
This section yields biholomorphic identifications $s : U \cong s(U)$ and $f : s(U) \cong U$.
Under these maps, the push-forward $h = f_{*}g$ restricted to $U$ can be identified with the restriction $g|_{s(U)}$ on $s(U)$.
This proves that $f_{*}g$ is smooth and Nakano (resp.~ Griffiths) semi-positive on $U$.
\end{proof}

We prove the constancy statement for metrics used above.
It is easier to discuss this in the dual setting.
Although we only need \ref{trivial} in the case where the metric $h$ is smooth in \ref{pushw}, we state it for possibly singular Hermitian metrics.
As mentioned earlier, we refer to \cite[\S\S 16--19]{HPS} for the relevant definitions and basic properties.

\begin{lem}\label{trivial} 
Let $X$ be a complex manifold, and let $E$ be a holomorphic vector bundle on $X$ of rank $r$. 
Suppose that $E$ admits a singular Hermitian metric $h$ with Griffiths semi-negative curvature. 
Let $Y \subset X$ be a reduced and connected compact complex subspace, and suppose $E$ is trivial on $Y$, i.e., $E|_{Y} = Y \times \BC^{r}$. 
Suppose that there exists a point $x \in Y$ such that $h_{x}$ is positive definite.
Then, $h|_{Y}$ is constant; i.e., there exists a positive definite Hermitian matrix $H_{0} (=h_{x})$ of size $r$ such that $h|_{Y} = Y \times H_{0}$. 
\end{lem}  

\begin{proof} 
We first assume that $Y=X$. 
The rank $1$ case is well known. 
The higher rank case can be reduced to the rank $1$ case as follows. 
We write $h = (h_{i\ol j})$. 
For each $i = 1, \ldots, r$, let $e_{i} \in H^{0}(X,E)$ be the constant section corresponding to the standard basis vector $(0, \ldots, 1, \ldots, 0)$ in $\BC^{r}$, where the $i$-th entry is $1$ and all other entries are $0$. 
Then, the squared norm $|e_{i}|_{h}^{2} = h_{i\ol i}$ is a plurisubharmonic (psh) function (possibly discontinuous) on $X$. 
Since $X$ is compact, $h_{i\ol i}$ must be constant. 

Next, we consider the section $s = e_{1} + e_{2} \in H^{0}(X,E)$. 
Its squared norm is given by $|s|_{h}^{2} = h_{1\ol 1} + h_{2\ol 2} + 2 \mathrm{Re}\, h_{1\ol 2}$. 
By the same argument, the psh function $|s|_{h}^{2}$ is constant, which implies that $\mathrm{Re}\, h_{1\ol 2}$ is constant. 
Similarly, for $t = e_{1} + \ai e_{2}$, we have $|t|_{h}^{2} = h_{1\ol 1} + h_{2\ol 2} + 2 \mathrm{Im}\, h_{1\ol 2}$. 
This shows that $\mathrm{Im}\, h_{1\ol 2}$ is also constant. 
By repeating this argument for all pairs of indices, we conclude that every entry $h_{i\ol j}$ is constant.  

Now, let us consider the general case. 
Let $Y = \bigcup_{i=1}^{k} Y_{i}$ be the irreducible decomposition of $Y$. 
We choose a smooth model $f_{i} : X_{i} \to Y_{i}$ for each $Y_{i}$. 
If $(\det h)|_{Y_{i}} \not\equiv 0$, the preceding argument implies that $f_{i}^{*}h = X_{i} \times H_{i}$ for some positive definite Hermitian matrix $H_{i}$. 
In this case, $\det h$ is nowhere zero on $Y_{i}$.

We may assume that $h_{x}$ is positive definite at some point $x \in Y_{1}$. 
Then, we initially have $h|_{Y_{1}} = Y_{1} \times h_{x}$. 
Since $Y$ is connected, we can find another component, say $Y_{2}$, such that $Y_{1} \cap Y_{2} \neq \emptyset$. 
It then follows that $(\det h)|_{Y_{2}} \not\equiv 0$, which implies $f_{2}^{*}h = X_{2} \times H_{2}$. 
By repeating this procedure, we can propagate this property to all components. 
Consequently, the desired assertion holds.
\end{proof}

Proposition \ref{pushw} yields the following slight generalization of \ref{mt}.
We outline the proof below.

\begin{thm}\label{mt-variant}
Let $X$ be a complex space, $Y$ a complex manifold, and $f:X\to Y$ a weakly semi-stable morphism with connected fibers.
Suppose that the toroidal structures of $X$ and $Y$ are strict, or that $X$ has at worst quotient singularities.
Let $M$ be a holomorphic vector bundle on $X$.
Suppose that for any relatively compact open subset $V \subset Y$,
there exists a modification $\mu : X_{V}' \to X_{V}=f^{-1}(V)$ from a complex manifold $X_{V}'$ such that
the composition $f_{V}' = f\circ \mu$ is locally K\"ahler and that $M_{V}' = \mu^{*}M$ is Nakano semi-positive.
Then, the direct image sheaf $E=f_{*}(K_{X/Y} \otimes M)$ is locally free, and the universal quotient line bundle $\OE1$ on $\BP(E)$ is metrically nef.
\end{thm}

\begin{proof}[Outline of proof]
Suppose $f:X\to Y$ is toroidal with respect to the toroidal structures $(U_{X}\subset X)$ and $(U_{Y}\subset Y)$.
For notational simplicity, $Y$ stands for a relatively compact open subset of the original $Y$.
By assumption, there exists a modification $\mu:X'\to X$ from a complex manifold $X'$ such that $f' = f\circ \mu : X' \to Y$ is locally K\"ahler and $M' := \mu^{*}M$ admits a Hermitian metric $h'$ with Nakano semi-positive curvature.
By \ref{pushw}, we have a $C^{0}$-Hermitian metric $h$ on $M$.

As $M'$ is Nakano semi-positive, thanks to the Grauert--Riemenschneider type vanishing (\cite[Theorem 5.2\,(iii)]{takegoshi}), we have $R^{i}\mu_{*}(K_{X'}\ot M') = 0$ for any $i>0$.
In particular, we have $R^{q}f'_{*}(K_{X'/Y}\ot M') = R^{q}f_{*}(K_{X/Y}\ot M)$ for any $q \ge 0$.
(Local freeness of these direct image sheaves follows by adapting the argument in \S 2.)
We are interested in 
$$
	E' := f'_{*}(K_{X'/Y}\ot M') =  f_{*}(K_{X/Y}\ot M) =: E. 
$$
Let $X_y$ (resp.~ $X'_y$) denote the fiber of $f$ (resp.~ $f'$) over a point $y \in Y$.
The modification $\mu: X' \to X$ induces a holomorphic map $\mu_y: X'_y \to X_y$.
Put
$$
	Y_{1} = \{ y \in U_{Y} \mid 
	\text{$f'$ is smooth over $y$ and $\mu_{y}:X'_{y} \to X_{y}$ is bimeromorphic} \}. 
$$
This is non-empty and analytic Zariski open in $Y$.

Let $g'$ be the canonical $L^{2}$-metric on $E'$ with respect to $h'$.
We know that $g'$ is a singular Hermitian metric on $E'$ with Griffiths semi-positive curvature.
We know $g'$ is $C^{\infty}$ and Griffiths (in fact Nakano) semi-positive on $E'|_{Y_{1}}$.
On the other hand, we can define a canonical $L^{2}$-metric $g$ on $E|_{U_{Y}}$ with respect to $h$.
Since $h$ is continuous, integration along the fibers of the smooth proper map $f|_{U_{X}}$ shows that $g$ is continuous on $U_{Y}$.
We now observe that 
$$
	g' = g \ \text{ on } Y_{1}.
$$
This is because, if $y \in Y_{1}$, $\mu_{y}:X'_{y}\to X_{y}$ is bimeromorphic and hence
$H^{0}(X'_{y}, K_{X'_{y}} \ot M') = H^{0}(X_{y}, K_{X_{y}} \ot M)$ via ${\mu_{y}}_{*}$. 
The proof of the local freeness of $E'$ over $Y_{1}$ also guarantees $E'_{y}=H^{0}(X'_{y}, K_{X'_{y}} \ot M')$ and $E_{y} = H^{0}(X_{y}, K_{X_{y}} \ot M)$.
The above identification is an isometry with respect to the natural $L^{2}$-pairings; equivalently, $g'_{y} = g_{y}$.
More concretely, if $s' \in H^{0}(X'_{y}, K_{X'_{y}} \ot M')$, we have $s = {\mu_{y}}_{*}s' \in H^{0}(X_{y}, K_{X_{y}} \ot M)$. 
Moreover, 
$$
	|s'|_{g'}^{2} = \int_{X'_{y}} h' {\rm i}_n s' \wed \ol {s'} 
	= \int_{X_{y}} h{\rm i}_n s \wed \ol {s} 
	= |s|_{g}^{2},
$$
where $n= \dim X - \dim Y$.

In particular, $g$ is $C^{\infty}$ on $Y_{1} (\subset U_{Y})$ (and $C^{0}$ on $U_{Y}$).
By the continuity of $g$ and uniqueness of psh extension applied to the dual metrics, the equality $g' = g$ extends to $U_{Y}$.
We henceforth denote the global metric $g'$ by $g$.
We see the logarithmic divergence of $g$ along $Y \sm U_{Y}$ in the same way as in the proof of \ref{mt}.
For this, we need to use the weakly semi-stable morphism $f:X\to Y$ and the fact that $h$ is $C^{0}$.
To see the logarithmic divergence of $g$, we only need the continuity and non-degeneracy of $h$ (the smoothness of $h$ is not needed).
Hence $g$ is a singular Hermitian metric on $E$ with Griffiths semi-positive curvature, which is $C^{0}$ on $U_{Y}$, and has at most logarithmic divergence along $Y \sm U_{Y}$.
This is enough to conclude the metric nefness of $\OE{1}$.
\end{proof}


\end{document}